\documentclass[twocolumn]{autart}

\makeatletter
\@ifundefined{cl@part}{\def\cl@part{}}{}%
\@ifundefined{cl@section}{\def\cl@section{}}{}%
\makeatother

\usepackage{balance}

\usepackage{pgfplots}
\pgfplotsset{compat=newest}
\usepackage{graphicx}
\usepackage{tikz}
\usepackage{overpic}
\usepackage{stanli}
\usetikzlibrary{arrows.meta,positioning,intersections,shapes,arrows}

\usepackage{amsmath,amssymb,amsfonts,mathtools,bm}

\usepackage[inline]{asymptote}

\usepackage{times}
\usepackage{xcolor}
\usepackage{enumitem}

\usepackage[numbers,sort&compress]{natbib}

\renewcommand*{\qed}{\hfill\ensuremath{\square}}

\usepackage[hidelinks]{hyperref}

\newtheorem{ex}[thm]{Example}

\DeclarePairedDelimiter{\norm}{\lVert}{\rVert}

\DeclarePairedDelimiterX{\setdef}[2]{\{}{\}}{#1\,\delimsize\vert\,\mathopen{} #2}
\DeclarePairedDelimiterX{\scprod}[2]{\langle}{\rangle}{#1,#2}
\DeclarePairedDelimiterX{\dualprod}[2]{\langle}{\rangle}{#1,#2}
\DeclarePairedDelimiterX{\sdprod}[2]{\llangle}{\rrangle}{#1,#2} 
\newcommand{\dd}{\mathrm{d}} \newcommand{\dx}[1][x]{\mathop{\dd#1}}

\newcommand{\Lb}{\mathcal{L}_{\mathrm{b}}} 
\newcommand{\Lp}[1]{\mathrm{L}^{#1}} 
\newcommand{\Wkp}[1]{\mathrm{W}^{#1}} 

\newcommand{\conC}{\mathrm{C}} 

\newcommand{\Hk}[1]{\mathrm{H}^{#1}} 

\newcommand{\adjunsymb}{\ast} 

\newcommand{\adjun}[1][1]{%
  \setcounter{i}{1}%
  \striche={\adjunsymb}%
  \loop%
  \ifnum\value{i}<#1%
  \striche=\expandafter{\the\expandafter\striche\adjunsymb}%
  \stepcounter{i}%
  \repeat%
  ^{\the\striche}%
}
\newcommand{\X}{\mathcal{X}}
\newcommand{\U}{\mathcal{U}}
\newcommand{\Y}{\mathcal{Y}}

\renewcommand{\d}{\mathrm{d}}

\newcommand{\loc}{\textup{loc}}
\renewcommand{\H}{\mathrm{H}}

\newcommand{\pmatsmall}[1]{\left(\begin{smallmatrix}#1 \end{smallmatrix}\right)}

\newcommand{\re}{\mathrm{Re}}

\newcommand{\Id}{I}

\newcommand{\ext}{\operatorname{ext}}

\renewcommand{\Re}{\operatorname{Re}}

\newcommand{\dom}{\operatorname{dom}}
\newcommand{\im}{\operatorname{im}}

\newcommand{\R}{\ensuremath{\mathbb R}}    
\newcommand{\C}{\ensuremath{\mathbb C}}    

\newcommand{\bvek}[2]
{
   \begin{bmatrix}
      #1\\
      #2
   \end{bmatrix}
}
\newcommand{\pvek}[2]
{
   \begin{pmatrix}
      #1\\
      #2
   \end{pmatrix}
}
\newcommand{\sbvek}[2]{\left[\begin{smallmatrix}#1\\#2\end{smallmatrix}\right]}
\newcommand{\spvek}[2]{\left(\begin{smallmatrix}#1\\#2\end{smallmatrix}\right)}
\newcommand{\sbmat}[4]{\left[\begin{smallmatrix}#1 & #2\\#3 & #4\end{smallmatrix}\right]}

\usepackage{graphicx}          

\begin{document}

\begin{frontmatter}

\title{Funnel control for passive infinite-dimensional systems\thanksref{footnoteinfo}} 

\thanks[footnoteinfo]{This paper was not presented at any IFAC 
meeting. Corresponding author A.~Hastir.\\
This work was supported by the Deutsche Forschungsgemeinschaft (DFG, German Research Foundation) -- Project numbers 362536361, \emph{Adaptive control of coupled rigid and flexible multibody systems with port-Hamiltonian structure} and 532208976, \emph{System theoretic properties of linear infinite-dimensional port-Hamiltonian systems}, and the Research Council of Finland grant number 349002.
}

\author[1]{Thavamani Govindaraj}\ead{thavamani.govindaraj@tu-ilmenau.de}  
\author[2]{Anthony Hastir}\ead{hastir@uni-wuppertal.de}  
\author[3]{Lassi Paunonen}\ead{lassi.paunonen@tuni.fi} 
\author[1]{Timo Reis}\ead{timo.reis@tu-ilmenau.de}  

\address[1]{Institut f\"ur Mathematik, Technische Universit\"at Ilmenau, Weimarer Stra\ss e 25, 98693 Ilmenau, Germany}  
\address[2]{Fachbereich Mathematik und Naturwissenschaften, Bergische Universit\"at Wuppertal, Gau\ss stra\ss e 20, 42119 Wuppertal, Germany}  
\address[3]{Faculty of Information Technology and Communication Sciences, Tampere University, Kalevantie 4, FI-33100 Tampere, Finland}


\begin{abstract}                          
We consider funnel control for linear infinite-dimensional systems that are impedance passive, meaning that they satisfy an energy balance in which the stored energy equals the squared norm of the state and the supplied power is the inner product of input and output. For the analysis we employ the system node approach, which offers a unified framework for infinite-dimensional systems with boundary and distributed control and observation. The resulting closed-loop dynamics are governed by a nonlinear evolution equation; we establish its solvability and hence the applicability of funnel control to this class. The applicability is illustrated by an Euler–Bernoulli beam, which is studied in two distinct scenarios: once with boundary control and once with distributed control.\end{abstract}

\end{frontmatter}

\section{Introduction}

We consider output regulation for input–output dynamical systems, where the objective is to design a controller such that the system output $y$ follows a given reference signal $y_{\mathrm{ref}}$. To this end, the controller generates the input $u$ based on the tracking error $e = y - y_{\mathrm{ref}}$, see Figure~\ref{fig:blockdiag}.

\tikzstyle{block} = [draw, thick, fill=white, rectangle, 
    minimum height=3em, minimum width=6em]
\tikzstyle{sum} = [draw, thick, fill=white, circle, node distance=1cm]
\tikzstyle{input} = [coordinate]
\tikzstyle{output} = [coordinate]
\tikzstyle{pinstyle} = [pin edge={to-,thick,black}]

\begin{figure}[ht]
    \centering
\begin{tikzpicture}[auto, node distance=2cm,>=latex',thick]

    \node [input, name=input] {};
    \node [sum, right of=input] (sum) {};
    \node [block, right of=sum] (controller) {Controller};
    \node [sum, right of=controller, xshift = 2em] (sum2) {}; 
    \node [block, right of=sum2, xshift = -1em] (system) {System};
    \node [output, right of=system] (output) {};

    \node at ([yshift=-4,xshift=-5]sum.west) {$-$};
    \node at ([yshift=-6,xshift=5]sum.south) {$+$};
    \node at ([yshift=-0.1em,xshift=-1em]sum2.south) {$+$};
    \node at ([yshift=6,xshift=-0.5em]sum2.north) {$+$};

    \draw [->] (input) -- node {$y_{\rm ref}$} (sum);
    \draw [->] (sum) -- node {$e$} (controller);
    \draw [->] (controller) -- (sum2);
    \draw [->] (sum2) -- node[name=u] {$u$} (system);
    \draw [->] (system) -- node [name=y] {$y$} (output);
    \draw [->] (y) -- ++(0,-2) -| (sum.south);

    \node [above of=sum2,yshift=-2em] (uext) {};
    \draw [->] (uext) -- node[right,yshift=0.4em] {$u_{\mathrm{ext}}$} (sum2);

\end{tikzpicture}
    \caption{Feedback interconnection of controller and system with external input $u_{\mathrm{ext}}$.}
    \label{fig:blockdiag}
\end{figure}
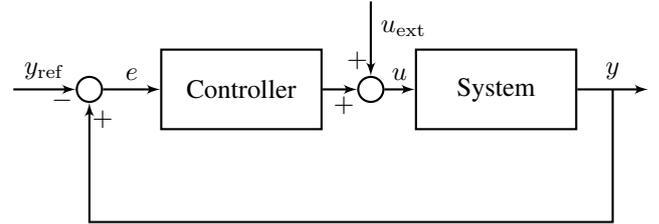

Funnel control is an adaptive error feedback strategy that constrains the tracking error $e$ to evolve within a prescribed performance funnel. Let $\U$ denote the space in which the outputs (and likewise the inputs) take values. Given a reference signal $y_{\mathrm{ref}}:\R_{\ge0}\to \U$, and assuming sufficient smoothness (to be specified later), the controller generates an input $u(t)$ such that 
\begin{equation*}
   \forall\,t\ge0:\quad \varphi(t)\|e(t)\| < 1
\end{equation*}
for some given positive $\varphi:\R_{\ge0}\to\R$ (the {\em funnel function}) with the property that $\varphi$ and $\varphi^{-1}$ are bounded.

A typical funnel control is given by
\begin{equation}
u(t) = -\frac{1}{1-\varphi(t)^2\|e(t)\|^2}e(t) + u_{\ext}(t),
\label{eq:FunnelLaw}
\end{equation}
where $u_{\ext}$ denotes an external input such as a disturbance, a feedforward control, or a combination of both. The precise role of this term will be discussed later in the paper. The essential part, however, is the feedback component represented by the first summand: if the error approaches the funnel boundary, i.e., $\|e(t)\|\approx 1/\varphi(t)$, then the denominator $1-\varphi(t)^2|e(t)|^2$ becomes small and the controller applies a large gain to drive the error back into the funnel. Conversely, if the error is close to zero, the control action is correspondingly small.

\begin{figure}[h]
\centering
\begin{tikzpicture}[scale=0.45]
\tikzset{>=latex}
  \filldraw[color=gray!25] plot[smooth] coordinates {(0.15,4.7)(0.7,2.9)(4,0.4)(6,1.5)(9.5,0.4)(10,0.333)(10.01,0.331)(10.041,0.3) (10.041,-0.3)(10.01,-0.331)(10,-0.333)(9.5,-0.4)(6,-1.5)(4,-0.4)(0.7,-2.9)(0.15,-4.7)};
  \draw[thick] plot[smooth] coordinates {(0.15,4.7)(0.7,2.9)(4,0.4)(6,1.5)(9.5,0.4)(10,0.333)(10.01,0.331)(10.041,0.3)};
  \draw[thick] plot[smooth] coordinates {(10.041,-0.3)(10.01,-0.331)(10,-0.333)(9.5,-0.4)(6,-1.5)(4,-0.4)(0.7,-2.9)(0.15,-4.7)};
  \draw[thick,fill=lightgray] (0,0) ellipse (0.4 and 5);
  \draw[thick] (0,0) ellipse (0.1 and 0.333);
  \draw[thick,fill=gray!25] (10.041,0) ellipse (0.1 and 0.333);
  \draw[thick] plot[smooth] coordinates {(0,2)(2,1.1)(4,-0.1)(6,-0.7)(9,0.25)(10,0.15)};
  \draw[thick,->] (-2,0)--(12,0) node[right,above]{\normalsize $t$};
  \draw[thick,dashed](0,0.333)--(10,0.333);
  \draw[thick,dashed](0,-0.333)--(10,-0.333);
  \node [black] at (0,2) {\textbullet};
  \draw[->,thick](4,-3)node[right]{\normalsize ${\displaystyle\lambda=
  \frac1{\sup_{t\ge0}\varphi(t)}}$}--(2.5,-0.4);
  \draw[->,thick](3,3)node[right]{\normalsize $(0,e(0))$}--(0.07,2.07);
  \draw[->,thick](9,3)node[right]{\normalsize $1/\varphi(t)$}--(7,1.4);
\end{tikzpicture}
\caption{Evolution of the error $e(t)$ inside a funnel with boundary $1/\varphi(t)$.}
\label{Fig:funnel}
\end{figure}
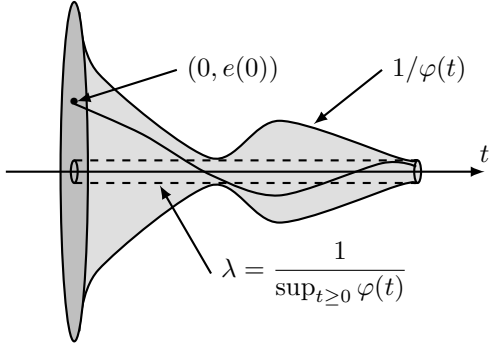

\medskip

Funnel control has been extensively studied in the finite-dimensional setting. After its introduction in \cite{IlcRyaSan2002}, numerous extensions followed, covering multi-input multi-output systems \cite{IlcRyaTre2005}, nonlinear systems with strict relative degree \cite{Ber2018}, and nonlinear systems subject to input constraints \cite{Ber2024}. Applications include chemical reactors \cite{IlcTre2004}, wind turbines \cite{Hac2014}, robotic manipulators \cite{HacKen2012}, electrical circuits \cite{BerRei2014}, and position control of mechanical systems \cite{HacHopIlcMueTre2013}. For surveys on the theory and applications of funnel control we refer to \cite{IlchmannRyan2008, BerIlcRya2021}.

\medskip

For infinite-dimensional systems, funnel control has so far only been studied in specific cases. Initial approaches transferred the proof strategy from the finite-dimensional setting, which led to the consideration of the Byrnes-Isidori form for infinite-dimensional systems \cite{IlchmannSeligTrunk2016}. Based on this, funnel control has been considered for a certain class of such systems, later extended to include systems with moderate nonlinearities \cite{HasWinDoc2023}. However, these approaches require highly restrictive assumptions on the control and observation operators, so that in particular boundary control is not covered. Nevertheless, there are practically relevant systems for which funnel control can be established within this framework, for instance a moving water tank \cite{Ber2020,Ber2022}.

Boundary control in the context of funnel control was first considered in \cite{ReiSel2015} for the heat equation. A more general class of boundary control problems was investigated in \cite{MarTimFel2021}, where, under a suitable passivity assumption, funnel control was established by employing the theory of maximally monotone operators and nonlinear evolution equations.
The underlying idea has also been applied to the nonlinear monodomain equation with boundary control \cite{BergerBreitenPucheReis2021}.
 
While boundary control for passive infinite-dimensional systems is by now well understood, the seemingly simpler situation of distributed control is not covered by the existing theory (in particular not by \cite{IlchmannSeligTrunk2016}, where even stronger smoothness assumptions on the control and observation operators are imposed). In our view, a general theory of funnel control for passive infinite-dimensional systems that treats both boundary and distributed control is highly desirable.

This is precisely the contribution of the present article. We develop a unifying approach to funnel control for passive infinite-dimensional linear systems that encompasses both boundary and distributed control. Our analysis builds on the system node framework, see e.g.~\cite[Sec.~4.7]{sta2005}, which provides a natural and flexible setting for treating boundary and distributed actuation and sensing within a single theory.
Systems within this framework are described by
\begin{equation}\label{eq:ODEnode}
  \spvek{\dot{x}(t)}{y(t)}
  \;=\;
  \sbvek{A\&B\\[-1mm]}{C\&D}\spvek{x(t)}{u(t)},
\end{equation}
where $A\&B:\X\times\U\supset\dom(A\&B)\to\X$ and $C\&D:\dom(A\&B)\to\Y$ are linear operators (with some further properties clarified in due course). 
The somewhat unusual use of the ``$\&$'' symbol emphasizes that the domain is in general not a Cartesian product of subspaces of $\X$ and $\U$. This becomes particularly relevant in boundary control problems. 
The assumptions of the theory are formulated so broadly that virtually any linear PDE with control and observation can be cast into this setting. 
In recent years a number of works have demonstrated the applicability of the framework to various complex examples \cite{PhilippReisSchaller2025,ReisSchaller2025LQOC,ReisSchaller2024OseenChapter,FarkasJacobReisSchmitz2023}.

\medskip

Section~\ref{sec:prelim} covers system-node basics; Section~\ref{sec:ass} states our and discusses our assumptions; Section~\ref{sec:main} proves feasibility of the closed-loop dynamics under funnel control. We conclude in Section~\ref{sec:ex} with an illustration of the results by means of an Euler-Bernoulli beam, considered in two scenarios: with boundary control and with distributed control.
\medskip

\subsection*{Notation}
The open right complex half-plane is denoted by $\C_+$.
All spaces considered in this article are complex Hilbert spaces. 
The norm in $\X$ and the inner product in $\X$ are denoted by $\norm{\cdot}_{\X}$ and $\scprod{\cdot}{\cdot}_{\X}$, respectively, or simply by $\norm{\cdot}$ and $\langle \cdot,\cdot \rangle$ when no ambiguity arises. The symbol $\Lb(\X,\Y)$ denotes the set of all bounded linear operators from $\X$ to $\Y$, and we abbreviate $\Lb(\X):=\Lb(\X,\X)$, and $I_\X\in\Lb(\X)$ is the identity (again, we skip the subindex, if it is clear from context).
For self-adjoint operators $P,Q\in\Lb(\X)$, we write $P\geq Q$, if $P-Q$ is a~nonnegative operator in the sense that $\scprod{x}{Px}_\X\geq0$ for all $x\in\X$. 

A linear operator $A:\X\supset\dom(A)\to\Y$ is said to be {\em closed} if $\dom(A)$ becomes a Hilbert space when endowed with the inner product
\[
    \scprod{x_1}{x_2}_{\dom(A)}
      := \scprod{x_1}{x_2}_{\X} + \scprod{Ax_1}{Ax_2}_{\Y},
\]
and {\em densely defined}, if $\dom(A)$ is a~dense subspace of $\X$.
The symbols $\im A$ and $\ker A$ denote the range and the kernel of $A$, respectively. 
For an operator $A:\X\supset\dom(A)\to\X$ we write $\rho(A)$ and $\sigma(A)$ for its resolvent set and its spectrum, respectively.

Throughout this work we follow the conventions for Lebesgue and Sobolev spaces as presented in the standard reference \cite{AdamsFournier2003}. 
Whenever function spaces with values in a Hilbert space $\X$ are considered, we indicate the target space explicitly by writing “$;\X$’’ after the domain. 
For instance, the space of $p$-integrable $\X$-valued functions on a domain $\Omega \subset \R^d$ is written as $\Lp{p}(\Omega;\X)$. 
Integration of such $\X$-valued functions is always understood in the sense of Bochner integration, cf.\ \cite{DiestelUhl1977}. 

\section{System nodes}\label{sec:prelim}

We start by collecting some preliminaries for systems of the form \eqref{eq:ODEnode}. 
The homogeneous part of the dynamics, that is, the case $u=0$, is governed by the so-called \emph{main operator}
\begin{align*}
  A\colon \X\supset\dom(A)&\to\X, \\
  \dom(A)&\coloneqq\setdef{x\in\X}{\spvek{x}{0}\in\dom(A\&B)},\\
  x&\mapsto A\&B\,\spvek{x}{0}.
\end{align*}
The framework of system nodes imposes structural conditions on $A\&B$ and $C\&D$ to ensure that the resulting dynamical system admits a consistent solution theory.  This includes the concept of strongly continuous semigroups (that is, the extension of the matrix exponential to operators), \cite{EngelNagel2000}. Note that our definition of system nodes differs slightly from that in \cite{sta2005}; however, the two are equivalent by \cite[Lem.~2.3]{PhilippReisSchaller2025}.

\begin{defn}[System node]\label{def:sysnode}
  A \emph{system node} on a triple $(\X,\U,\Y)$ of Hilbert spaces is a linear operator $S=\sbvek{A\&B}{C\&D}$ with
  \begin{align*}
    A\&B:\X\times\U\supset\dom(A\&B)&\to\X,\\
    C\&D:\X\times\U\supset\dom(C\&D)&\to\Y,
  \end{align*}
  that satisfies:
  \begin{enumerate}[label=(\alph{*})]
    \item $A\&B$ is closed.
    \item $C\&D\in\Lb(\dom(A\&B),\Y)$.
    \item For each $u\in\U$ there exists $x\in\X$ such that $\spvek{x}{u}\in\dom(A\&B)$.
    \item The associated main operator $A$ generates a strongly continuous semigroup $\mathbb{T}(\cdot)\colon\R_{\ge0}\to\Lb(\X)$.
  \end{enumerate}
\end{defn}

\begin{defn}[Classical and generalized trajectories]\label{def:traj}
Let $\sbvek{A\&B}{C\&D}$ be a system node on $(\X,\U,\Y)$ and $T>0$.  
\begin{itemize}
  \item A \emph{classical trajectory} of \eqref{eq:ODEnode} on $[0,T]$ is a triple
  \begin{multline*}
(x(\cdot),u(\cdot),y(\cdot))\\\in\conC^{1}([0,T];\X)\times \conC([0,T];\U)\times \conC([0,T];\Y),
  \end{multline*}
  that satisfies \eqref{eq:ODEnode} for every $t\in[0,T]$.
  \item A \emph{generalized trajectory} on $[0,T]$ is a limit of classical trajectories in the space
  \[
  \conC([0,T];\X)\times \Lp{2}([0,T];\U)\times \Lp{2}([0,T];\Y).
  \]
\end{itemize}
By extension, $(x,u,y)$ is called a classical or generalized trajectory on $\R_{\ge0}$ if its restriction to each finite interval $[0,T]$ is of the corresponding type.  
\end{defn}

\begin{prop}[{\cite[Lem.~4.7.8]{sta2005}}]\label{prop:solex}
Let $S=\sbvek{A\&B}{C\&D}$ be a system node on $(\X,\U,\Y)$, let $T>0$, and suppose that
$x_0\in\X$ and $u\in\Wkp{2,1}([0,T];\U)$ with $\spvek{x_0}{u(0)}\in\dom(A\&B)$.
Then there exist $x(\cdot)\in\conC^1([0,T];\X)$ with $x(0)=x_0$ and $y(\cdot)\in\conC([0,T];\Y)$ such that $(x,u,y)$ is a classical trajectory of \eqref{eq:ODEnode} on $[0,T]$.
\end{prop}

\begin{rem}[Further properties of system nodes]\label{rem:nodes}
Let $\sbvek{A\&B}{C\&D}$ be a system node on $(\X,\U,\Y)$.  
\begin{enumerate}[label=(\alph{*})]
\item\label{rem:nodes1} $A\&B$ extends to an operator \[[A_{-1}\,B]\in\Lb(\X_{-1}\times\U,\X),\]
where $\X_{-1}$ is the completion of $\X$ with respect to $\norm{x}_{\X_{-1}}:=\norm{(\alpha I-A)^{-1}x}$ for some $\alpha\in\rho(A)$. The topology of $\X_{-1}$ does not depend on $\alpha$ \cite[Prop.~2.10.2]{TucsnakWeiss2009}. Then $\dom(A\&B)$ fulfills
\begin{equation}
  \dom(A\&B)=\setdef{\spvek{x}{u}\in\X\times\U}{A_{-1}x+Bu\in\X},\label{eq:ABform}
\end{equation}
see \cite[Def.~4.7.2 and Lem.~4.7.3]{sta2005}. In what follows, we denote by the same letter $A$ the operator $A_{-1}$.

\item\label{rem:nodes5} We can introduce $C \in \Lb(\dom(A),\Y)$ with
\[
   Cx \;=\; C\&D \spvek{x}{0}.
\]
\end{enumerate}
\end{rem}

We now introduce the notion of the transfer function of a system node. 
To this end, recall from \cite[Lem.~4.7.3]{sta2005} that 
\[
   \sbvek{(\lambda I - A_{-1})^{-1}B}{I}\in\Lb(\U;\dom(A\&B)).
\]

\begin{defn}[Transfer function]\label{def:TranfFct_S}
Let $S = \sbvek{A\& B}{C\& D}$ be a~system node on $(\X,\U,\Y)$. The \emph{transfer function $P$ associated to $S$} is
\[\begin{aligned}
P\colon &&\rho(A)\to &\,\mathcal{L}(\U,\Y),\\&& \lambda\mapsto&\,C\& D \sbvek{(\lambda\Id-A_{-1})^{-1}B}\Id.
\end{aligned}\]
\end{defn}
In this article, we focus on systems which are passive in the following sense. 
\begin{defn}[Impedance passive systems]\label{def:passive}
Let $S=\sbvek{A\&B}{C\&D}$ be a system node on $(\X,\U,\U)$ (i.e., the input and output spaces are equal). The system \eqref{eq:ODEnode} is called {\em impedance passive}, if for all $T>0$, all classical (and thus also all generalized) trajectories $(x(\cdot),u(\cdot),y(\cdot))$ on $[0,T]$, we have
\[\norm{x(T)}_\X^2-\norm{x(0)}_\X^2\leq 2\Re\int_0^T\scprod{u(t)}{y(t)}_\U\dx[t].\]
\end{defn}
This property has been studied in detail in \cite{Staffans2002Passive,ArovStaffans2005KYP,Staffans2002Personal}, see also \cite{Reis2025dissipation}. In particular, it follows from \cite[Thm.~4.1]{Staffans2002Personal} that impedance passivity is equivalent to
    \begin{multline}
\forall\,\spvek{x_0}{u_0} \in\dom(A\&B):\\
    \Re\scprod {A \& B \spvek{x_0}{u_0}}{x_0}_\X  \leq \Re\scprod{C \& D{\spvek{x_0}{u_0}}}{u_0}_\U.\label{eq:KYP}
    \end{multline}
Note that \eqref{eq:KYP} means that the operator $\sbvek{\phantom{-}A\& B}{-C\& D}$ is dissipative (and thus also maximal dissipative by \cite[Thm.~4.1]{Staffans2002Personal}). Let us give some further comments on this system class.
\begin{rem}[Impedance passive systems]\label{rem:nodes_ImpPass}
Let $\sbvek{A\&B}{C\&D}$ be a system node on $(\X,\U,\U)$ with transfer fucntion $P$. \begin{enumerate}[label=(\alph{*})]
\item\label{rem:nodes_ImpPass1} By setting $u\equiv 0$, we obtain that the semigroup $\mathbb{T}$ 
associated to an impedance passive system is contractive, i.e., $\norm{\mathbb{T}(t)}_{\Lb(\X)}\leq1$ for all $t\ge0$.
Then \cite[Prop.~2.3.1]{TucsnakWeiss2009} yields that $\C_+\subset\rho(A)$.
As~a consequence, the transfer function $P$ of an  impedance passive system is defined on a~superset of $\C_+$.
\item\label{rem:nodes_ImpPass2} If \eqref{eq:ODEnode} is impedance passive, then
\[P(\lambda)+P(\lambda)^*\geq 0\quad\forall\,\lambda\in\C_+.\]
In the finite-dimensional setting, this property is referred to as {\em positive realness} \cite{AndersonVongpanitlerd1973}. 
\item\label{rem:nodes_ImpPass3} An important subclass of impedance passive systems is given by 
{\em port-Hamiltonian system nodes} \cite{PhilippReisSchaller2025}. 
These systems can be written in the form
\begin{equation*}
  \spvek{\dot{x}(t)}{y(t)}
  \;=\;
  \sbvek{F\&G\\[-1mm]}{K\&L}\spvek{Hx(t)}{u(t)}.
\end{equation*} 
The operator $H\in\Lb(\X,\X^*)$ is assumed to be self-dual, positive, surjective, and 
\[
   \sbvek{F\&G}{K\&L}:\X^*\times \U\supset\dom(F\&G)\to \X\times\U
\] 
is dissipative, where $\X^*$ denotes the anti-dual of $\X$, that is, the space of bounded conjugate-linear functionals on $\X$. Such systems are impedance passive when $\X$ is endowed with the norm 
$\norm{x} := \scprod{x}{Hx}_{\X}$. This class covers a wide range of PDE systems, including Maxwell’s equations, advection–diffusion equations, Euler–Bernoulli and Timoshenko beams, the wave equation, and the Oseen equation \cite{ReisSchaller2024OseenChapter,FarkasJacobReisSchmitz2023,ReisSchaller2025LQOC,PhilippReisSchaller2025}.
\end{enumerate}
\end{rem}

\section{Our assumptions}\label{sec:ass}
Here we present and discuss the 
assumptions on the system class for which funnel control 
will be treated in Section \ref{sec:main}. 
\begin{assum}[System class]\label{assum:SystemClass}\
    \begin{enumerate}[label=(\alph{*})]
\item $\X$ and $\U$ are complex Hilbert spaces.
\item $S=\sbvek{A\&B}{C\&D}$ is a~system node on $(\X,\U,\U)$, such that the system \eqref{eq:ODEnode} is impedance passive.
\item   
    The operator $C:\dom(A)\to\U$ (see Remark~\ref{rem:nodes}\,\ref{rem:nodes5}) is surjective.
\item\label{assum:Invert_P} There exist some $\lambda,c>0$, such that the transfer function associated to $S$ satisfies
\[P(\lambda)+P(\lambda)^*\geq c\, I_\U.\]
\end{enumerate}
\end{assum}
We now turn to a discussion of our assumptions. We begin by characterizing the surjectivity of $C$.

\begin{prop}\label{prop:Csurj}
Let $\sbvek{A\&B}{C\&D}$ be a~system node on $(\X,\U,\Y)$. Then the following are equivalent:
\begin{enumerate}[label=(\alph{*})]
\item\label{prop:Csurj1} $C$ is surjective.
\item\label{prop:Csurj2} For all $u_0\in\U$, $y_0\in\Y$, there exists some $x_0\in\X$, such that
\[\spvek{x_0}{u_0}\in\dom(A\&B)\;\mbox{and}\;y_0=C\&D\spvek{x_0}{u_0}.\]
\item\label{prop:Csurj3} There exists  $Q\in \Lb(\Y,\dom(A))$ with $CQ=I_\Y$. 
\end{enumerate}
\end{prop}
\begin{pf}
The equivalence between \ref{prop:Csurj1} and \ref{prop:Csurj3} is an immediate consequence of the results in \cite{Erdelyi1972GeneralizedInverse} on generalized inverses of operators in Hilbert spaces.\\ 
    \ref{prop:Csurj1}$\Rightarrow$
    \ref{prop:Csurj2}: Assume that $u_0\in\U$, $y_0\in\Y$. Then, by definition of a~system node, there exists some $x_1\in\X$ with $\spvek{x_1}{u_0}\in\dom(A\&B)$. Surjectivity of $C$ yields that there exists some $x_2\in\dom(A)$ with
    $Cx_2=y_0-C\&D\spvek{x_1}{u_0}$.
Then a~simple calculation yields that $x_0:=x_1+x_2$ has the desired properties.\\
    \ref{prop:Csurj2}$\Rightarrow$
    \ref{prop:Csurj1}: Let $y_0\in\Y$. Then, by setting $u_0:=0$, \ref{prop:Csurj2} yields that there exists some $x_0$, such that $\spvek{x_0}0\in\dom(A\&B)$ and $C\&D\spvek{x_0}0=y_0$. This means that $y_0=Cx_0$, and the result is proven.\qed
\end{pf}
\begin{rem}[System class]\label{rem:sysclass}\
    \begin{enumerate}[label=(\alph{*})]
\item\label{rem:sysclass1} Our assumptions do not involve {\em well-posedness}. This notion refers to the continuous dependence of 
$y \in \Lp{p}(\R_{\ge0};\Y)$ on the initial value $x(0)=x_0 \in \X$ and the input 
$u \in \Lp{p}(\R_{\ge0};\U)$, $p \in [1,\infty]$. 
This omission is deliberate: well-posedness is not a natural property in the context of impedance passive systems (and would therefore exclude many relevant examples); further,
 verifying well-posedness in concrete cases is typically an~involved task.
\item\label{rem:sysclass3} Abstract boundary control systems are of the form
\begin{equation}
    \dot{x}(t)=\mathcal{A}x(t),\;\;
    u(t)=\mathcal{K}x(t),\;\;
    y(t)=\mathcal{L}x(t),\label{eq:bcs}
\end{equation}
where $\mathcal{A}:\X\supset\dom(\mathcal{A})\to\X$ is closed, $\mathcal{K}\in\Lb(\dom(\mathcal{A}),\U)$, $\mathcal{L}\in\Lb(\dom(\mathcal{A}),\Y)$, and the restriction of $\mathcal{A}$ to $\ker \mathcal{K}$ is the generator of a~strongly continuous semigroup. 
These can be turned into systems governed by system nodes \cite[Sec.~5.2]{sta2005}. Statement \ref{prop:Csurj2} in Proposition~\ref{prop:Csurj} (which is equivalent to surjectivity of $C$) is, in case of boundary control systems, equivalent to surjectivity of $\sbvek{\mathcal{K}}{\mathcal{L}}\in\Lb(\dom(\mathcal{A}),\U\times\Y)$. The latter condition has been also imposed in \cite{MarTimFel2021} in the context of funnel control of boundary control systems.
\item\label{rem:sysclass5} By using \cite[Cor.~4.5]{Logemann2020Spectral}, Assumption~\ref{assum:SystemClass}~\ref{assum:Invert_P} implies that the transfer function has the even stronger property that, for all $\lambda\in\C_{+}$, there exists some $c_\lambda>0$ with
\[
   P(\lambda)+P(\lambda)^* \;\geq\; c_\lambda I_\U.
\]
In particular, $P(\lambda)$ admits a bounded inverse for all $\lambda\in\C_+$, and it is straightforward to verify that the pointwise inverse is again positive real. 
\item\label{rem:sysclass6} If a system node $S=\sbvek{A\&B}{C\&D}$ on $(\X,\U,\U)$ has a positive real transfer function and, in addition, the system \eqref{eq:ODEnode} is both approximately controllable and approximately observable, then \cite[Thm.~5.1.1]{ArovStaffans2005KYP} implies that the state space can be renormed (not necessarily equivalently) such that an impedance passive is obtained.
\end{enumerate}
\end{rem}

We now state our assumptions on the reference signal $y_{\rm ref}$, on the funnel function~$\varphi$, on the initial state $x(0)=x_0$, and on the external input $u_{\ext}$ in the funnel control law
\eqref{eq:FunnelLaw}.
\begin{assum}\label{assum:Funnel_Ref_Signal}
\noindent\textbf{\upshape (Funnel, reference signal, initial state, and external input)}
\begin{enumerate}[label=(\alph{*})]
  \item The reference signal satisfies $y_{\mathrm{ref}} \in \Wkp{2,\infty}(\R_{\ge 0};\U)$.
  \item The funnel function satisfies $\varphi \in \Wkp{2,\infty}(\R_{\ge 0};\R)$, and there exists $\mu > 0$ such that $\varphi(t) \ge \mu$ for all $t>0$.
   \item The initial state $x_0\in\X$, the external input $u_{\ext}\in\Wkp{2,\infty}(\R_{\ge0};\U)$ and the funnel function $\varphi$ are such that there exists $e_0\in\U$ such that $\varphi(0)\Vert e_0\Vert < 1$, and
\begin{subequations}\label{eq:OutputInit}
  \begin{equation}
      \pvek{x_0}{u_{\ext}(0)-\frac{e_0}{1-\varphi(0)^2\Vert e_0\Vert^2}}\in \dom(A\&B)\label{eq:funnelinit}
  \end{equation}
with
\begin{equation}
      e_0 = C\&D\pvek{x_0}{u_{\ext}(0)-\frac{e_0}{1-\varphi(0)^2\Vert e_0\Vert^2}}-y_{\mathrm{ref}}(0).
\label{eq:OutputInit2}
  \end{equation}
\end{subequations}
\end{enumerate}
\end{assum}
\begin{rem}[Funnel assumptions]\label{rem:funnelass}\
\begin{enumerate}[label=(\alph{*})]
\item\label{rem:funnelass2} The boundedness of the funnel function implies the existence of a band around the reference signal that is entirely contained within the funnel, see Figure \ref{Fig:funnel}.
Moreover, strict boundedness of $\varphi$ from below means that the funnel region is bounded.
\item We will see in Theorem \ref{thm:Main_Thm_Exist_Uniq} that $\varphi(0)\|e_0\|<1$ simply means that the initial output is inside the funnel.
\item\label{rem:funnelass3} Condition~\eqref{eq:OutputInit} expresses that the controller is properly initialized at $t=0$. If the plant arises from the boundary control system \eqref{eq:bcs}, this amounts to a boundary consistency requirement: the boundary trace $\mathcal{K}x_0$ of the initial state must coincide with the value prescribed by the controller at $t=0$. An analogous compatibility condition already appears in \cite{MarTimFel2021} for funnel control of boundary control problems. Note also that \eqref{eq:OutputInit} imposes a constraint on the admissible initial state. In PDE settings, this typically entails sufficient spatial regularity so that the boundary trace $\mathcal{K}x_0$ is well defined and matches the controller’s initial signal. Note that, if the input operator is bounded (that is, $B$ in Remark~\ref{rem:nodes}\,\ref{rem:nodes1} is in $\Lb(\U,\X)$), then \eqref{eq:funnelinit} reduces to $x_0\in\dom(A)$.
\item\label{rem:funnelass4} As noted in the introduction, $u_{\ext}$ may represent an exogenous disturbance. It can also be a precomputed feedforward signal that guides the plant along a nominal trajectory, thereby reducing the burden on the feedback controller. 
This combination has proved particularly effective for mechanical systems; see, e.g., \cite{Druecker2025FunnelFFW,Berger2021Tracking}. If the external input $u_{\ext}$ comprises a precomputed feedforward term, an exogenous disturbance, or any linear combination of the two, one cannot, in general, expect \eqref{eq:OutputInit} to hold. To remedy this, the controller can be augmented to generate an additional signal that compensates the resulting mismatch at $t=0$. Rather than using only the standard funnel-feedback term, the controller may produce the signal
\begin{multline*}
-\frac{1}{1-\varphi(t)^2\|e(t)\|^2}e(t)\\
+ p(t)\left(\frac{e(0)}{1-\varphi(0)^2\|e(0)\|^2}-u_{\ext}(0)+u_0\right),
\end{multline*}
where $p\in \Wkp{2,\infty}(\R_{\ge 0})$ has compact support (hence it is eventually inactive), $p(0)=1$, and $u_0\in\U$ with $\spvek{x_0}{u_0}\in\dom(A\&B)$. This requires that $x_0$ is in
\[
\mathcal{V}
:= \setdef*{x_0\in\X}{\exists\,u_0\in\U \text{ s.t.\ }\spvek{x_0}{u_0}\in\dom(A\&B)},
\]
which, as mentioned above, typically corresponds to sufficient spatial regularity of the initial value in the PDE case. With this modification, the theory applies upon replacing $u_{\ext}$ by
\[
u_{\ext}(\cdot) + p(\cdot)\Bigl(\frac{e(0)}{1-\varphi(0)^2\|e(0)\|^2}-u_{\ext}(0)+u_0\Bigr).
\]
\end{enumerate}
\end{rem}

The next assumption is not required for feasibility of funnel control, but it is used as an additional condition to ensure that all signals in the closed-loop system are globally bounded.

\begin{assum}\label{ass:passshift}
    $\sbvek{A\&B}{C\&D}$ is a~system node on $(\X,\U,\U)$, and there exists some $\alpha>0$, such that the system 
    \begin{equation*}
    \spvek{\dot{x}(t)}{y(t)}
  \;=\;
  \left(\sbvek{A\&B\\[-1mm]}{C\&D}+\frac\alpha2\sbvek{I\phantom{\&}0\\[-1mm]}{0\phantom{\&}0}\right)\spvek{x(t)}{u(t)}
      \end{equation*}
is impedance passive.
\end{assum}

\begin{rem}
    Let $S=\sbvek{A\&B}{C\&D}$ be a system node on $(\X,\U,\U)$. 
The result 
\cite[Thm.~4.1]{Staffans2002Personal}
implies that Assumption~\ref{ass:passshift} is equivalent to 
    \begin{multline}
\forall\,\spvek{x_0}{u_0} \in\dom(A\&B):\\
    \Re\scprod {A \& B \spvek{x_0}{u_0}}{x_0}_\X  \\\leq -\alpha \norm{x_0}_\X^2+\Re\scprod{C \& D{\spvek{x_0}{u_0}}}{u_0}_\U.
    \label{eq:DissipIneq}
    \end{multline}
\end{rem}

\section{Main Results}\label{sec:main}

We start by presenting the statement of our main result on the applicability of funnel control. We observe that the interconnection between \eqref{eq:ODEnode} and the funnel control law \eqref{eq:FunnelLaw} leads formally to the nonlinear system
\begin{subequations}
\begin{align}
    \dot{x}(t) &= A\&B\left(\begin{smallmatrix}
        x(t)\\
        u_{\ext}(t)-\frac{e(t)}{1-\varphi(t)^2\Vert e(t)\Vert^2}
    \end{smallmatrix}\right),\label{eq:ClosedLoopODE}\\
    y(t) &= C\& D\left(\begin{smallmatrix}x(t)\\u_{\ext}(t)-\frac{e(t)}{1-\varphi(t)^2\Vert e(t)\Vert^2}\end{smallmatrix}\right),\label{eq:Output_ClosedLoop}\\
    \varphi(t)\Vert e(t)\Vert &< 1,\qquad e(t) = y(t) - y_\mathrm{ref}(t).\label{eq:FunnelPerf_ClosedLoop}
\end{align}
\end{subequations}
Our main result below shows that this system has a solution in an appropriate sense and that the funnel control achieves output tracking. 
\begin{thm}\label{thm:Main_Thm_Exist_Uniq}
Suppose that Assumption \ref{assum:SystemClass} holds and that $y_{\mathrm{ref}}, \varphi, u_{\ext}$ and $x_0$ satisfy Assumption \ref{assum:Funnel_Ref_Signal}. Then, the funnel control \eqref{eq:FunnelLaw} achieves output tracking when applied to \eqref{eq:ODEnode}. More precisely, there exist unique $x\in \Wkp{1,\infty}_\loc(\R_{\ge0};\X)$ and $y\in\Lp{\infty}_\loc(\R_{\ge0};\U)$ such that
\begin{itemize}
    \item for every $t\geq 0$, we have 
    \begin{equation*}
    Ax(t) + B\left(u_{\ext}(t)-\frac{e(t)}{1-\varphi(t)^2\Vert e(t)\Vert^2}\right)\in \X
    \end{equation*}
    and $y(t)$ is the unique solution of 
    \begin{equation*}
    y(t) = C\& D\left(\begin{smallmatrix}x(t)\\ u_{\ext}(t)-\frac{e(t)}{1-\varphi(t)^2\Vert e(t)\Vert^2}\end{smallmatrix}\right)
    \end{equation*}
    with $\varphi(t)\Vert e(t)\Vert < 1$, where $e(t) = y(t)-y_{\mathrm{ref}}(t)$;
    \item \eqref{eq:ClosedLoopODE}--\eqref{eq:FunnelPerf_ClosedLoop} hold for a.e. $t\geq 0$.
\end{itemize}
Moreover, the funnel control law $u(t) = u_{\ext}(t) - \frac{e(t)}{1-\varphi(t)^2\Vert e(t)\Vert^2}$ satisfies $u\in\Lp{\infty}_\loc(\R_{\ge0};\U)$.
\end{thm}

Before turning to the proof of our main result, we first collect several auxiliary lemmas that will be needed in the proof. 
   Let $B_1(0) := \setdef{u \in \U}{\|u\|_\U < 1}$.
In what follows, we make use of the function $\phi:B_1(0)\subset \U\to\U$ defined by \begin{equation}
    \phi(w) = \frac{w}{1-\Vert w\Vert^2}.
    \label{eq:NonLin_Phi}
\end{equation}
This function is monotone in the sense that
\begin{multline}
    \forall \,w_1,w_2\in B_1(0)\text{ with }w_1\neq w_2:\\ \re\scprod{w_1-w_2}{\phi(w_1)-\phi(w_2)}>0.\label{eq:phimon}
\end{multline}
The non-strict inequality has been proven in  \cite[Lem.~5.4]{MarTimFel2021}. The strict inequality follows by a~careful inspection of that proof by further invoking that $r\mapsto\frac{r}{1-r^2}$ is strictly monotonically increasing on the open interval $(-1,1)$.
We begin with a lemma on the solvability of the nonlinear equation 
\begin{equation}
   w = r - P\phi(w),
   \label{eq:ImpEq}
\end{equation}
for $w \in B_1(0)$, where $r \in \U$ is fixed and $P\in \Lb(\U)$ is coercive. 

\begin{lem}\label{lem:ExistUniq_ImpEq}
Assume that $P\in\Lb(\U)$ with $P+P^*\geq cI_\U$ for some $c>0$.
Further, let $\phi: B_1(0)\subset \U\to\U$ be defined in \eqref{eq:NonLin_Phi}. Then, for every $r\in \U$, the equation \eqref{eq:ImpEq} has a unique solution $w\in B_1(0)$. Moreover, the resulting solution map $\U\to B_1(0)$, $r\mapsto w(r)$, is globally Lipschitz continuous.
\end{lem}
\begin{pf}
Our assumption implies that $P$ has a bounded inverse and $P^{-1}+(P^{-1})^*\ge c_0I_\U$ for some $c_0>0$.
Using this and~\eqref{eq:phimon} we have that
\begin{align*}
\MoveEqLeft c_0\norm{w_1-w_2}_\U^2\le
\Re\scprod{w_1-w_2}{P^{-1}(w_1-w_2)}_\U\\
&\quad+\Re\scprod{w_1-w_2}{\phi(w_1)-\phi(w_2)}_\U\\
&=\Re\scprod{w_1-w_2}{P^{-1}(w_1-w_2)+\phi(w_1)-\phi(w_2)}_\U, 
\end{align*}
which further implies 
\begin{multline}
   \forall\,w_1,w_2:\quad \norm{w_1+P\phi(w_1) -\big(w_2+P\phi(w_2)\big)}_\U\\\geq c_0\|P^{-1}\|^{-1}\,\norm{w_1-w_2}_\U.\label{eq:invLip}
\end{multline}
This estimate in particular shows that the solution of~\eqref{eq:ImpEq} must be unique and that the solution map is globally Lip\-schitz continuous. 
Finally, to prove existence of a solution, let $r\in\U$ and consider the map $\psi:B_1(0)\to \U$, $\psi(w):=P^{-1}w+\phi(w)$. 
Equation~\eqref{eq:ImpEq} is equivalent to $P^{-1}r=\psi(w)$. 
Since $P$ is invertible with $P^{-1}+(P^{-1})^*\ge c_0I_\U$,
the operator $\xi I+P^{-1}$ is invertible for all $\xi>0$, and we define the continuous function
\[
   p(\xi):=\|(\xi I+P^{-1})^{-1}P^{-1}r\|^2-1+\xi^{-1}, \quad \xi\ge1.
\]
Clearly $p(1)\ge0$ and $\lim_{\xi\to\infty}p(\xi)=-1$, so by the intermediate value theorem there exists some $\xi^*\geq1$ with $p(\xi^*)=0$. Set $w:=(\xi^*I+P^{-1})^{-1}P^{-1}r$. Then
   $\|w\|^2=1-{\xi^*}^{-1}<1$ and
   $\phi(w)=\xi^*w=\xi^*(\xi^*I+P^{-1})^{-1}P^{-1} r$.
Hence
\begin{multline*}
   \psi(w)=P^{-1}w+\phi(w)\\=(\xi^*I+P^{-1})(\xi^*I+P^{-1})^{-1}P^{-1}r=P^{-1} r.
\end{multline*}
Thus, $w\in B_1(0)$ satisfies $P^{-1}r=P^{-1}w+\phi(w)$, which is equivalent to \eqref{eq:ImpEq}.\qed
\end{pf}

The next lemma will play an important role in the proof of our main result.

\begin{lem}\label{lem:Regularity_z}
    Let $S = \left[\begin{smallmatrix}A\& B\\ C\&D\end{smallmatrix}\right]$ be a system node with transfer function $P$ for which Assumption \ref{assum:SystemClass} is satisfied. Consider the nonlinear operator $A_\phi: \dom(A_\phi)\subset \X\to\X$ defined as
     \begin{subequations}
\begin{align*}
    &A_\phi z = Az - B\phi(w(z)),
    \\
    &\dom(A_\phi) = \left\{z\in\X\,\vert\, \exists w\in\U,Az - B\phi(w)\in\X,\right.\nonumber\\
    &\left.\hspace{2cm}w = C\& D\left(\begin{smallmatrix}z\\-\phi(w)\end{smallmatrix}\right), \Vert w\Vert<1\right\}
\end{align*}  
\end{subequations}
where $w(z)$ is the element corresponding to $z$ in the definition of $\dom(A_\phi)$. Then $A_\phi$ is a single-valued $m$-dissipative operator. Moreover, for any $T>0$, $z_0\in \dom(A_\phi)$, $f\in\Wkp{1,\infty}([0,T];\X)$, and $\omega\in\Wkp{1,\infty}([0,T])$, there exists a unique $z\in\Wkp{1,\infty}([0,T];\X)$ such that $z(0) = z_0$, $z(t)\in\dom(A_\phi)$ for all $t\in [0,T]$ and
    \begin{equation}
        \dot{z}(t) = A_\phi(z(t)) + \omega(t)z(t) + f(t)        \label{eq:ODE_z_Proof}
        \end{equation}
        holds for a.e. $t\in [0,T]$. In addition, $\dot{z}(t)$ and $A_\phi(z(t))$ are continuous except at a countable number of values in $[0,T]$.
\end{lem}

\begin{pf}
Singled-valuedness and $m-$dissipativity of $A_\phi$ follows, by invoking \eqref{eq:phimon}, analogously as in~\cite[Proof of Thm.~3.9]{HasPau2025}. It then follows from~\cite[Thm.~4.20]{Miy92book} that the operator $A_\phi$ is the generator of a nonlinear contraction semigroup $\mathbb{T}^\phi(t)$ on $\overline{\dom(A_\phi)}$. Further, $\phi(0) = 0$ implies that $0\in \dom(A_\phi)$ and $A_\phi(0) = 0$. Since $f\in\Wkp{1,\infty}([0,\infty);\X)$, it satisfies in particular $f\in\Wkp{1,\infty}([0,T];\X)$ for every $T>0$. The claim then follows by \cite[Lem.~5.7]{MarTimFel2021}.\qed
\end{pf}

\begin{rem}\label{rem:Single-Valued}
The operator $A_\phi$ in Lemma~\ref{lem:Regularity_z} is single-valued if and only if the element $w$  in the definition of $\dom(A_\phi)$ is unique. 
The "if" part of the statement is clear. To verify the converse implication, 
suppose that $A_\phi $ is single-valued and assume that $w_1,w_2\in B_1(0)$ are two elements corresponding to $z\in \dom(A_\phi)$. Then $B(\phi(w_1)-\phi(w_2))=A_\phi z-A_\phi z=0\in X$. Thus $(0,\phi(w_1)-\phi(w_2))\in \dom(A\&B)$ and the passivity of the system implies
\begin{align*}
    0 
    &= \re \ \scprod{A\&B \left(\begin{smallmatrix}
        0\\ \phi(w_1)-\phi(w_2)
    \end{smallmatrix}\right)}{0}\\
    &\le \re \ \scprod{C\&D \left(\begin{smallmatrix}
        z-z\\ \phi(w_1)-\phi(w_2)
    \end{smallmatrix}\right)}{\phi(w_1)-\phi(w_2)}\\
    &= \re \ \scprod{w_2-w_1}{\phi(w_1)-\phi(w_2)}.
\end{align*}
However, \eqref{eq:phimon} now implies that we must have $w_1=w_2$.
\end{rem}

\begin{pf*}{Proof of Theorem~\textup{\ref{thm:Main_Thm_Exist_Uniq}}.}
By Proposition \ref{prop:Csurj}\,\ref{prop:Csurj3} we can choose $Q\in\Lb(\U,\dom(A))$ such that $CQ = I_\U$. Moreover, define $\omega$ and $f$ by $\omega(t) := \frac{\dot{\varphi}(t)}{\varphi(t)}$ and
\begin{align}
    f(t) &:= \varphi(t)\Big(AQ y_\mathrm{ref}(t) + \big(\lambda R_\lambda B- AQP(\lambda)\big)u_{\ext}(t)\nonumber\\
    &- Q\dot{y}_\mathrm{ref}(t) -\big(R_\lambda B-QP(\lambda)\big)\dot{u}_{\ext}(t)\Big)\label{eq:f}
    \end{align}
    where $R_\lambda = (\lambda I-A)^{-1}\in\Lb(\X_{-1},\X)$, $\lambda >0$. Now observe that thanks to $Q\in\Lb(\U,\dom(A))$ and $A\in\Lb(\dom(A),\X)$, we have $AQ\in\Lb(U,\X)$. Moreover, $Q P(\lambda)-R_\lambda B\in\Lb(\U,\X)$ and $\lambda R_\lambda B - AQP(\lambda)\in\Lb(\U,\X)$. 
    The properties $y_\mathrm{ref}\in\Wkp{2,\infty}(\R_{\ge0};\U)$, $\varphi\in\Wkp{2,\infty}(\R_{\ge0})$, $u_{\ext}\in\Wkp{2,\infty}(\R_{\ge0};\U)$ and $\omega\in\Wkp{1,\infty}(\R_{\ge0})$ therefore imply that $f\in\Wkp{1,\infty}(\R_{\ge0};\X)$. 
    Now we fix an arbitrary $T>0$, define 
    \begin{equation}
        z_0 = \varphi(0)\left(x_0 - Qy_{\mathrm{ref}}(0) - (R_\lambda B - QP(\lambda))u_{\ext}(0)\right)
        \label{eq:z0_def}
\end{equation}
and denote by $e_0\in\U$ an element that satisfies \eqref{eq:OutputInit2} in Assumption \ref{assum:Funnel_Ref_Signal}. We have $z_0\in\dom(A_\phi)$ with $w(z_0) := \varphi(0)e_0$ being the element `$w$' in the definition of $\dom(A_\phi)$. This element is unique according to Remark \ref{rem:Single-Valued}. Lemma \ref{lem:Regularity_z} implies that there exists a unique $z\in\Wkp{1,\infty}([0,T];\X)$ such that $z(0) = z_0$, $z(t)\in\dom(A_\phi)$ for all $t\geq 0$ and \eqref{eq:ODE_z_Proof}
holds for a.e. $t\in [0,T]$. In addition, $\dot{z}(t)$ and $A_\phi(z(t))$ are continuous except at a countable number of values in $[0,T]$. Expanding \eqref{eq:ODE_z_Proof} implies
\begin{equation}
\dot{z}(t) = Az(t) - B\phi(w(z(t))) + \omega(t)z(t) + f(t),
\label{eq:ODE_z_Inter}
\end{equation}
where $w(z(t))$ is the unique solution to
\begin{align}
w(z(t)) = C\& D\left(\begin{smallmatrix}
z(t)\\ -\phi(w(z(t)))
\end{smallmatrix}\right)
\label{eq:ImpEqu_wz}
\end{align}
with $\Vert w(z(t))\Vert < 1$. Uniqueness of $w(z(t))$ is guaranteed because $A_\phi$ is a single-valued operator, see Remark \ref{rem:Single-Valued}. Now for all $t\in[0,T]$
we define $x(t)$ and $y(t)$ by
\begin{align}
x(t) &:= \varphi(t)^{-1}z(t) + Q y_{\mathrm{ref}}(t) + (R_\lambda B - QP(\lambda))u_{\ext}(t)\label{eq:x_based_on_z}\\
y(t) &:= \varphi(t)^{-1}w(z(t)) + y_\mathrm{ref}(t) 
\label{eq:y_based_on_z}
\end{align}
 and denote $e(t) = y(t) - y_{\mathrm{ref}}(t)$. 
For every $t\in [0,T]$ we have (using $AR_\lambda B = \lambda R_\lambda B - B$ and $\varphi(t)e(t) = w(z(t))$)
\begin{align}
\MoveEqLeft[.5] Ax(t) + B(u_{\ext}(t) - \varphi(t)^{-1}\phi(\varphi(t)e(t)))\nonumber\\
&= \varphi(t)^{-1}Az(t) + AQy_\mathrm{ref}(t) + AR_\lambda B u_{\ext}(t)\nonumber\\
&\quad - AQP(\lambda)u_{\ext}(t)+ Bu_{\ext}(t) - \varphi(t)^{-1}B\phi(\varphi(t)e(t))\nonumber\\
&= \varphi(t)^{-1}\left(Az(t) - B\phi(w(z(t)))\right) + AQy_\mathrm{ref}(t)\nonumber\\[1ex]
&\quad+(\lambda R_\lambda B - AQP(\lambda))u_{\ext}(t).
\label{eq:AxBu_in_X}
\end{align}
Thus $Ax(t) + B(u_{\ext}(t) - \varphi(t)^{-1}\phi(\varphi(t)e(t))\in \X$ for all $t\in[0,T]$.
Now we compute $C\& D\left(\begin{smallmatrix}x(t)\\u_{\ext} - \varphi(t)^{-1}\phi(\varphi(t)e(t))\end{smallmatrix}\right)$. 
Using the properties of system nodes, \eqref{eq:y_based_on_z} and $CQ = I_\U$ we obtain
{\allowdisplaybreaks\begin{align}
\MoveEqLeft[.5] C\& D\left(\begin{smallmatrix}x(t)\\u_{\ext} - \varphi(t)^{-1}\phi(\varphi(t)e(t))\end{smallmatrix}\right)\nonumber\\
&=C\& D\left(\begin{smallmatrix}\varphi(t)^{-1}z(t) + Qy_{\mathrm{ref}}(t) + (R_\lambda B - QP(\lambda))u_{\ext}(t)\\u_{\ext} - \varphi(t)^{-1}\phi(\varphi(t)e(t))\end{smallmatrix}\right)\nonumber\\
&=\varphi(t)^{-1}C\& D\left(\begin{smallmatrix}z(t)\\-\phi(\varphi(t)e(t))\end{smallmatrix}\right) + CQ y_{\mathrm{ref}}(t)\nonumber\\
&\quad +C\& D\left(\begin{smallmatrix}
(R_\lambda B - QP(\lambda))u_{\ext}(t)\\
u_{\ext}(t)
\end{smallmatrix}\right)
\nonumber\\
&= \varphi(t)^{-1}C\& D\left(\begin{smallmatrix}
z(t)\\-\phi(w(z(t)))
\end{smallmatrix}\right) + y_\mathrm{ref}(t)\nonumber\\
&= \varphi(t)^{-1}w(z(t)) + y_{\mathrm{ref}}(t)
= y(t).\label{eq:y_from_w}
\end{align}%
}%
The estimate $\Vert w(z(t))\Vert < 1$ also implies $\varphi(t) \Vert e(t)\Vert<1$. This together with \eqref{eq:y_from_w} shows that \eqref{eq:Output_ClosedLoop}--\eqref{eq:FunnelPerf_ClosedLoop} hold for all $t\in [0,T]$. Now we focus on the regularity of $x$. According to the regularity of $z$ and the regularity of $\varphi$, $y_\mathrm{ref}$, $u_{\ext}$ in Assumption \ref{assum:Funnel_Ref_Signal} we get that $x\in\Wkp{1,\infty}([0,T];\X)$ since $Q\in\Lb(\U,\dom(A))$ and $R_\lambda B - QP(\lambda)\in\Lb(\U,\X)$. 
In the next step of the proof, we will show that $x$ satisfies \eqref{eq:ClosedLoopODE}. 
For $t\in [0,T]$ a direct computation using the formula for $f$ in \eqref{eq:f} shows that
\begin{align*}
    \dot{x}(t) &= -\frac{\dot{\varphi}(t)}{\varphi(t)}\varphi(t)^{-1}z(t) + \varphi(t)^{-1}\dot{z}(t)\\
    &\quad+ Q\dot{y}_{\mathrm{ref}}(t) + (R_\lambda B - QP(\lambda))\dot{u}_{\ext}(t)\\
    &=-\omega(t)\varphi(t)^{-1}z(t) + A\varphi(t)^{-1}z(t)\\
    &\quad- B\varphi(t)^{-1}\phi(w(z(t)))+\omega(t)\varphi(t)^{-1}z(t) \\
    &\quad+\! \varphi(t)^{-1}f(t)\!+\! Q\dot{y}_{\mathrm{ref}}(t) \!+\! (R_\lambda B \!-\! QP(\lambda))\dot{u}_{\ext}(t)\\
    &= A\varphi(t)^{-1}z(t) - \varphi(t)^{-1}B\phi(w(z(t)))\\
    &\quad +AQy_{\mathrm{ref}}(t) + (\lambda R_\lambda B - AQP(\lambda))u_{\ext}(t).
\end{align*}
The identity \eqref{eq:AxBu_in_X} finally implies that
\begin{align*}
    \dot{x}(t) = Ax(t) + B(u_{\ext}(t) - \varphi(t)^{-1}\phi(\varphi(t)e(t))),
\end{align*}
which is exactly \eqref{eq:ClosedLoopODE}.
To prove the uniqueness of $x$ and $y$, 
assume that these functions are as in the claim.
We can then define $z(t) $ and $w(z(t))$
by inverting the formulas \eqref{eq:x_based_on_z} and \eqref{eq:y_based_on_z}. Then $z\in\Wkp{1,\infty}([0,T];\X)$ and $\|w(z(t))\|<1$ for all $t\in [0,T]$, and computations similar to those in the first part of the proof show that  
\eqref{eq:ODE_z_Inter} and \eqref{eq:ImpEqu_wz} hold and $z(0)=z_0$, where $z_0$ is determined by \eqref{eq:z0_def}. The solutions of these two equations are unique by Lemma \ref{lem:Regularity_z}, and this implies the uniqueness of $x $ and $y$.

Now we study the regularity of $y$ defined in \eqref{eq:y_based_on_z}. 
The identity \eqref{eq:ODE_z_Inter} implies that 
\begin{align*}
z(t) &= R_\lambda(\lambda z(t) - \dot{z}(t) + \omega(t)z(t) + f(t))\nonumber\\
&\quad- R_\lambda B\phi(w(z(t))).
\end{align*}
Substituting this expression into
 \eqref{eq:ImpEqu_wz} yields
\begin{align}
    w(z(t)) &= C\& D\left(\begin{smallmatrix}z(t)\\-\phi(w(z(t)))\end{smallmatrix}\right)= r(t) - P(\lambda)\phi(w(z(t)))\label{eq:ImpEqua_Ter}
\end{align}
where $r(t) := CR_\lambda((\lambda + \omega(t)) z(t) - \dot{z}(t) + f(t))$. Since $CR_\lambda\in\Lb(\X,\U)$ and $(\lambda + \omega) z - \dot{z} + f\in\Lp{\infty}([0,T];\X)$, we have $r\in\Lp{\infty}([0,T];\U)$. Moreover, according to Lemma \ref{lem:ExistUniq_ImpEq}, we have that for every $t\in[0,T]$ there exists a~unique solution $w(z(t))\in B_1(0)$ of \eqref{eq:ImpEqua_Ter}. Lemma \ref{lem:ExistUniq_ImpEq} also implies that the map $(I+P(\lambda)\phi)^{-1}:\U\to B_1(0)$ is globally Lipschitz continuous. Since $(I+P(\lambda)\phi)^{-1}(0) = 0$, this in particular means that $\Vert w(z(t))\Vert\leq l_\lambda\Vert r(t)\Vert$ for some constant $l_\lambda>0$ depending only on $\lambda$. Consequently, $y$ defined in \eqref{eq:y_based_on_z} satisfies $y\in\Lp{\infty}([0,T];\U)$. We end the proof by studying the
regularity of the funnel control law $u(t) = u_{\ext}(t) - \varphi(t)^{-1}\phi(\varphi(t)e(t))$. Observe that $u(t) = u_{\ext}(t) - \varphi(t)^{-1}\phi(w(z(t)))$. According to $w(z)\in\Lp{\infty}([0,T];\U)$, $\phi(w(z)) = P(\lambda)^{-1}(r-w(z))\in \Lp{\infty}([0,T];\U)$ and $u_{\ext}\in\Wkp{2,\infty}(\R_{\ge0};\U)$, we conclude that $u\in\Lp{\infty}([0,T];\U)$. Since $T>0$ was arbitrary, the proof is complete.\qed
\end{pf*}

In Theorem \ref{thm:Main_Thm_Exist_Uniq}, we showed that funnel control is feasible, and in particular, that the solution to \eqref{eq:ClosedLoopODE} exists globally. However, there are even simple (even finite-dimensional) cases where the state $x$ and the control $u$ are not bounded.
\begin{ex}
\noindent\textbf{\upshape (Passive system with unbounded funnel control)}\\
    Consider the system 
    \begin{align*}
        \dot{x}(t)&=-x(t)+u(t),\\y(t)&=-x(t)+u(t).
    \end{align*}
The system is passive due to $\sbmat{\phantom{-}A}{\phantom{-}B}{{-}C}{-D}=\sbmat{-1}{\phantom{-}1}{\phantom{-}1}{-1}\leq0$. Let us consider the constant reference signal $y_{\rm ref}\equiv 1$, together with a constant funnel (i.e., a tube) of diameter 1, that is, $\varphi \equiv 2$. Our subsequent considerations are structural in nature; they are based neither on numerical computations nor on an explicit solution of the closed-loop differential equation. Theorem~\ref{thm:Main_Thm_Exist_Uniq} guarantees the existence of a global solution of the closed-loop system, with the additional property that $y(t)\in(1/2,3/2)$ for all $t\ge 0$. Then, for all $t\geq0$, $\tfrac12\leq y(t)=-x(t)+u(t)=\dot{x}(t)$, and thus
\[x(t)\geq x(0)+\tfrac{t}2.\]
Hence, $x$ is unbounded. To also see that $u$ is unbounded, we note that
\[u(t)=x(t)+y(t)\geq x(t)+\tfrac12\geq x(0)+\tfrac{t}2+\tfrac12.\]
\end{ex}
The following theorem guarantees the boundeness of both the state and the input of the closed-loop system \eqref{eq:ClosedLoopODE}--\eqref{eq:FunnelPerf_ClosedLoop} under an additional assumption.
\begin{thm}
Suppose that Assumptions \ref{assum:SystemClass} and \ref{ass:passshift} hold and that $y_{\mathrm{ref}}$, $\varphi$, $u_{\ext}$ and $x_0$ satisfy Assumption \ref{assum:Funnel_Ref_Signal}. Then the solution of the closed-loop system \eqref{eq:ClosedLoopODE}--\eqref{eq:FunnelPerf_ClosedLoop} satisfies $x\in\Wkp{1,\infty}(\R_{\ge0};\X)$. Moreover, the funnel control law $u(t) = u_{\ext}(t)-\varphi(t)^{-1}\phi(\varphi(t)e(t))$ satisfies $u\in\Lp{\infty}(\R_{\ge0};\U)$.
\end{thm}

\begin{pf}
In what follows, $Q\in\Lb(\U,\dom(A))$ is such that $CQ = I_\U$. We start by considering \eqref{eq:ODE_z_Inter} which governs the evolution of $z = \varphi\left(x - Q y_\mathrm{ref}-(R_\lambda B  - QP(\lambda))u_{\ext}\right)\in\Wkp{1,\infty}_\loc(\R_{\ge0};\X)$, where $f\in\Wkp{1,\infty}(\R_{\ge0};\X)$ is given in \eqref{eq:f} and $\omega = \frac{\dot{\varphi}}{\varphi}\in\Wkp{1,\infty}(\R_{\ge0})$. According to the relation between $z$ and $x$ in \eqref{eq:x_based_on_z}, it is easy to see that $x\in\Wkp{1,\infty}(\R_{\ge0};\X)$ if and only if $z\in\Wkp{1,\infty}(\R_{\ge0};\X)$. We start by showing that $z\in\Lp{\infty}(\R_{\ge0};\X)$. For this, observe that
\begin{align*}
    &\frac{1}{2}\frac{\d}{\d t}\Vert z(t)\Vert^2_\X = \re\langle \dot{z}(t),z(t)\rangle_\X \\
    &= \re\langle A\& B\pmatsmall{z(t)\\ -\phi(w(z(t)))},z(t)\rangle_\X + \omega(t)\Vert z(t)\Vert^2_\X\\
    & + \re\langle f(t),z(t)\rangle_\X\\
    &\leq \re\langle C\& D\pmatsmall{z(t)\\-\phi(w(z(t)))},-\phi(w(z(t)))\rangle_\U\\
    &+ (\omega(t)-\alpha)\Vert z(t)\Vert_\X^2 + \Vert f(t)\Vert_\X\Vert z(t)\Vert_\X\\
    &= -\re\langle w(z(t)),\phi(w(z(t)))\rangle_\U + (\omega(t)-\alpha)\Vert z(t)\Vert_\X^2\\
    &+ \Vert f(t)\Vert_\X\Vert z(t)\Vert_\X\\
    &\leq (\omega(t)-\alpha)\Vert z(t)\Vert_\X^2 + \Vert f(t)\Vert_\X\Vert z(t)\Vert_\X.
\end{align*}
According to \cite[Lem.~5.8]{MarTimFel2021} and $\omega = \frac{\dot{\varphi}}{\varphi} = \frac{\d}{\d t}\log(\varphi)$, we can estimate $\Vert z(t)\Vert_\X$ in the following way
\begin{align*}
    \Vert z(t)\Vert_\X&\leq \varphi(t)e^{-\alpha t}\left(\frac{\Vert z(0)\Vert_\X}{\varphi(0)}+\int_0^t\frac{e^{\alpha s}}{\varphi(s)}\Vert f(s)\Vert_\X\d s\right).
\end{align*}
Using that $\varphi, \varphi^{-1}\in\Lp{\infty}(\R_{\ge0})$ and $f\in\Wkp{1,\infty}(\R_{\ge0};\X)$, we get that
\begin{align*}
    \Vert z(t)\Vert_\X\leq \frac{\Vert z(0)\Vert_\X}{\varphi(0)}\Vert\varphi\Vert_\infty + \frac{\Vert \varphi\Vert_\infty}{\alpha}\Vert\varphi^{-1}\Vert_\infty\Vert f\Vert_\infty,
\end{align*}
whence $z\in\Lp{\infty}(\R_{\ge0};\X)$. Now we will show that $\dot{z}\in\Lp{\infty}(\R_{\ge0};\X)$. Let $h>0$ and denote $\phi_{t+h} := -\phi(w(z(t+h)))$, $\phi_t := -\phi(w(z(t)))$. Then, we have
\begin{align*}
    &\tfrac{1}{2}\tfrac{\d}{\d t}\Vert z(t+h) - z(t)\Vert^2_\X\\
    &= \re\langle \dot{z}(t+h) - \dot{z}(t),z(t+h)-z(t)\rangle_\X\\
    &= \re\langle A\& B\left(\begin{smallmatrix}z(t+h)-z(t)\\\phi_{t+h}-\phi_t\end{smallmatrix}\right), z(t+h)-z(t)\rangle_\X\\
    & + \re\langle \omega(t+h)z(t+h)-w(t)z(t),z(t+h)-z(t)\rangle_\X\\
    &+\re\langle f(t+h)-f(t),z(t+h)-z(t)\rangle_\X.
\end{align*}
Thanks to the dissipation inequality \eqref{eq:DissipIneq} we have that
{\allowdisplaybreaks\begin{align*}
    &\tfrac{1}{2}\tfrac{\d}{\d t}\Vert z(t+h)-z(t)\Vert^2_\X\\
    &\leq \re\langle C\& D\left(\begin{smallmatrix}z(t+h)-z(t)\\\phi_{t+h}-\phi_t\end{smallmatrix}\right),\phi_{t+h}-\phi_t\rangle_\U\\
    &+(\omega(t+h)-\alpha)\Vert z(t+h)-z(t)\Vert^2_\X\\
    &+ \vert \omega(t+h)-\omega(t)\vert\Vert z(t)\Vert_\X\Vert z(t+h)-z(t)\Vert_X\\
    & + \Vert f(t+h)-f(t)\Vert_\X\Vert z(t+h)-z(t)\Vert_\X.
    \end{align*}}
Because $z\in\Lp{\infty}(\R_{\ge0};\X)$ and 
\begin{align*}
&\re\langle C\& D\left(\begin{smallmatrix}z(t+h)-z(t)\\\phi_{t+h}-\phi_t\end{smallmatrix}\right),\phi_{t+h}-\phi_t\rangle_\U\\
&= \re\langle w(z(t+h))-w(z(t)),\phi_{t+h}-\phi_t\rangle_\U\leq 0,
\end{align*}
we have
\begin{align*}
    &\tfrac{1}{2}\tfrac{\d}{\d t}\Vert z(t+h)-z(t)\Vert^2_\X\\
    &\leq (\omega(t+h)-\alpha)\Vert z(t+h)-z(t)\Vert^2_\X\\
    &+ \vert \omega(t+h)-\omega(t)\vert\Vert z\Vert_\infty\Vert z(t+h)-z(t)\Vert_X\\
    & + \Vert f(t+h)-f(t)\Vert_\X\Vert z(t+h)-z(t)\Vert_\X.
\end{align*}
Applying \cite[Lem.~5.8]{MarTimFel2021} to the previous inequality implies that
\begin{align*}
&\Vert z(t+h)-z(t)\Vert_\X\leq \Vert z(h)-z(0)\Vert_\X e^{\int_0^t\omega (s+h)-\alpha\d s}\\
&+\int_0^t e^{\int_s^t \omega(r+h)-\alpha)\d r}\Big(\vert\omega(s+h)-\omega(s)\vert\Vert z\Vert_\infty\\
&\hspace{3cm}+ \Vert f(s+h)-f(s)\Vert_\X\Big)\d s\\
&= \Vert z(h)-z(0)\Vert_\X\frac{\varphi(t+h)}{\varphi(h)}e^{-\alpha t}\\
&+\varphi(t+h)\int_0^t \frac{e^{-\alpha(t-s))}}{\varphi(s+h)}\Big(\vert\omega(s+h)-\omega(s)\vert\Vert z\Vert_\infty\\
&\hspace{3cm}+ \Vert f(s+h)-f(s)\Vert_\X\Big)\d s.
\end{align*}
Dividing the previous estimate by $h>0$ and letting $h\to 0$ implies that
\begin{align*}
    &\Vert\dot{z}(t)\Vert_\X\leq \frac{\Vert \dot{z}(0)\Vert_\X}{\varphi(0)}\varphi(t)e^{-\alpha t}\\
    &+ \varphi(t)\int_0^t\frac{e^{-\alpha(t-s)}}{\varphi(s)}\Big( \vert \dot{w}(s)\vert\Vert z\Vert_\infty + \Vert\dot{f}(s)\Vert_\X\Big)\d s\\
    &\leq \frac{\Vert A_\phi(z_0) + \omega(0)z_0 + f(0)\Vert_\X}{\varphi(0)}\Vert\varphi\Vert_\infty e^{-\alpha t}\\
    &+\frac{\Vert\varphi\Vert_\infty}{\Vert\varphi^{-1}\Vert_\infty}\Big( \vert\dot{\omega}\vert_\infty\Vert z\Vert_\infty + \Vert\dot{f}\Vert_\infty\Big)\frac{1-e^{-\alpha t}}{\alpha}.
\end{align*}
From this, we deduce that $\dot{z}\in\Lp{\infty}(\R_{\ge0};\X)$. It remains to show that $u(t) = u_{\ext}(t) - \varphi(t)^{-1}\phi(\varphi(t)e(t))$ satisfies $u\in\Lp{\infty}(\R_{\ge0};\U)$. Thanks to the regularity of $u_{\ext}$, this reduces to showing that $\varphi^{-1}\phi(\varphi e)\in\Lp{\infty}(\R_{\ge0};\U)$. According to the notations used in the proof of Theorem~\ref{thm:Main_Thm_Exist_Uniq}, we have $\varphi(t)^{-1}\phi(\varphi(t)e(t)) = \varphi(t)^{-1}\phi(w(z(t)))$. Moreover, according to Lemma \ref{lem:ExistUniq_ImpEq}, we have that for all $t\geq 0$, $w(z(t))$ satisfies $w(z(t))\in B_1(0)$, which implies that $w(z)\in\Lp{\infty}(\R_{\ge0};\U)$. Now we consider \eqref{eq:ImpEqua_Ter}, where $r(t) = CR_\lambda((\lambda + \omega(t)) z(t) - \dot{z}(t) + f(t))$.  From this equation, we get that $   \phi(w(z(t))) = P(\lambda)^{-1}\left(r(t)-w(z(t))\right)$, which, according to $w(z)\in\Lp{\infty}(\R_{\ge0};\U)$ and $r\in\Lp{\infty}(\R_{\ge0};\U)$, implies that $u = -\varphi^{-1}\phi(w(z)) + u_{\ext}\in\Lp{\infty}(\R_{\ge0};\U)$.\qed
\end{pf}

\section{Example: An Euler-Bernoulli beam}\label{sec:ex}
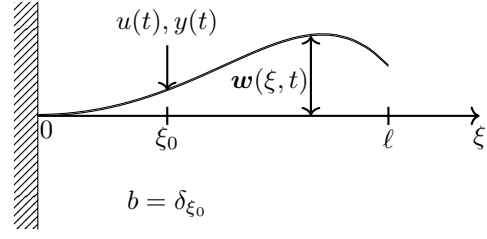
\begin{figure}
\centering
\begin{tikzpicture}[scale=0.6]

\tikzstyle{ground}=[fill,pattern=north east lines,draw=none,minimum width=0.75cm,minimum height=0.3cm]

\node (wall) [ground, rotate=-90, minimum width=3cm,yshift=1.05cm] {};
\draw (wall.north east) -- (wall.north west);

\def\xmin{2}
\def\xmax{9.7}
\def\w{8}
\def\nsamples{20}

\draw[thick,->] (\xmin,0)--(\xmax+2,0) node[right,below]{\normalsize $\xi$};
  \draw[variable=\x,samples=\nsamples,smooth,domain=\xmin:\xmax]
    plot(\x,{0.07*(\x-2)^2*cos(0.022*(\x-2)^2 r)});
    \draw[variable=\x,samples=\nsamples,smooth,domain=\xmin:\xmax]
    plot(\x,{0.07*(\x-2)^2*cos(0.022*(\x-2)^2 r)+0.03});
    
\draw[thick,<->] (\w,0)--(\w,{0.07*(\w-2)^2*cos(0.022*(\w-2)^2 r)}) {};
\node[thick] at (\w-0.9,{0.07*(\w-2)^2*cos(0.022*(\w-2)^2 r)-1}){\normalsize ${\bm w}(\xi,t)$};
\node[thick] at (\xmin+0.2,-0.3) {\normalsize $0$};
\draw[thick,-] (\xmax,-0.2) -- (\xmax,0.2) {};
\node[thick] at (\xmax,-0.5) {\normalsize $\ell$};

\draw[thick,-] (\xmax/2,-0.2) -- (\xmax/2,0.2) {};
\node[thick] at (\xmax/2,-0.5) {\normalsize $\xi_0$};
\draw[thick,->] (\xmax/2,{0.07*(\xmax/2-2)^2*cos(0.022*(\xmax/2-2)^2 r)+1}) -- (\xmax/2,{0.07*(\xmax/2-2)^2*cos(0.022*(\xmax/2-2)^2 r)});
\node[thick] at (\xmax/2,{0.07*(\xmax/2-2)^2*cos(0.022*(\xmax/2-2)^2 r)+1.5}) {\normalsize $u(t),y(t)$};

\node[thick] at (\xmax/2,{0.07*(\xmax/2-2)^2*cos(0.022*(\xmax/2-2)^2 r)-2.5}) {\normalsize $b = \delta_{\xi_0}$};

\end{tikzpicture}
\caption{A Beam with pointwise control and observation, i.e., $b = \delta_{\xi_0}$.}
	\label{fig:beam_bnd}
\end{figure}

\begin{figure}
\centering
\begin{tikzpicture}[scale=0.6]

\tikzstyle{ground}=[fill,pattern=north east lines,draw=none,minimum width=0.75cm,minimum height=0.3cm]

\node (wall) [ground, rotate=-90, minimum width=3cm,yshift=1.05cm] {};
\draw (wall.north east) -- (wall.north west);

\def\xmin{2}
\def\xmax{9.7}
\def\w{8}
\def\x{\xmax/2-0.35}
\def\xt{\xmax/2+0.35}
\def\nsamples{20}

\draw[thick,->] (\xmin,0)--(\xmax+2,0) node[right,below]{\normalsize $\xi$};
  \draw[variable=\x,samples=\nsamples,smooth,domain=\xmin:\xmax]
    plot(\x,{0.07*(\x-2)^2*cos(0.022*(\x-2)^2 r)});
    \draw[variable=\x,samples=\nsamples,smooth,domain=\xmin:\xmax]
    plot(\x,{0.07*(\x-2)^2*cos(0.022*(\x-2)^2 r)+0.03});
    
\draw[thick,<->] (\w,0)--(\w,{0.07*(\w-2)^2*cos(0.022*(\w-2)^2 r)}) {};
\node[thick] at (\w-0.9,{0.07*(\w-2)^2*cos(0.022*(\w-2)^2 r)-1}){\normalsize ${\bm w}(\xi,t)$};
\node[thick] at (\xmin+0.2,-0.3) {\normalsize $0$};
\draw[thick,-] (\xmax,-0.2) -- (\xmax,0.2) {};
\node[thick] at (\xmax,-0.5) {\normalsize $\ell$};

\draw[thick,->] (\x,{0.07*(\x-2)^2*cos(0.022*(\x-2)^2 r)+1}) -- (\x,{0.07*(\x-2)^2*cos(0.022*(\x-2)^2 r)});
\draw[thick,->] (\xt,{0.07*(\xt-2)^2*cos(0.022*(\xt-2)^2 r)+1}) -- (\xt,{0.07*(\xt-2)^2*cos(0.022*(\xt-2)^2 r)});
\draw[thick,->] (\xmax/2,{0.07*(\xmax/2-2)^2*cos(0.022*(\xmax/2-2)^2 r)+1}) -- (\xmax/2,{0.07*(\xmax/2-2)^2*cos(0.022*(\xmax/2-2)^2 r)});
\node[thick] at (\xmax/2,{0.07*(\xmax/2-2)^2*cos(0.022*(\xmax/2-2)^2 r)+1.5}) {\normalsize $u(t),y(t)$};
\node[thick] at (\xmax/2,{0.07*(\xmax/2-2)^2*cos(0.022*(\xmax/2-2)^2 r)-2.5}) {\normalsize $b\in\Lp{2}([0,\ell])$};
  \end{tikzpicture}
  \caption{A Beam with distributed control and observation, i.e. $b\in\Lp{2}([0,\ell])$.}
	\label{fig:beam_distr}
\end{figure}
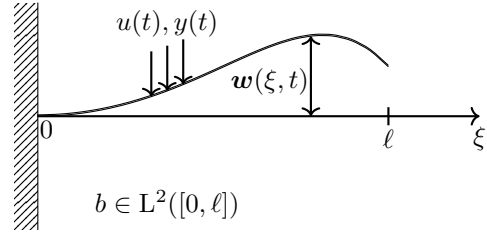

In this section we illustrate the main results from Section~\ref{sec:main} using an Euler--Bernoulli beam that is clamped at its left end. The equation for the beam is (at least formally) given by
\begin{equation*}
\begin{aligned}
\rho(\xi)\tfrac{\partial^2}{\partial t^2}{\bm w}(\xi,t)&=-\tfrac{\partial^2}{\partial\xi^2}\big(EI(\xi)\tfrac{\partial^2{\bm w}}{\partial\xi^2}(\xi,t)\big)+b(\xi)u(t),\\
y(t)&=\int_0^\ell b(\xi)
\tfrac{\partial}{\partial t}{\bm w}(\xi,t)\dx[\xi],\quad t\ge0,\,\xi\in [0,\ell],
\end{aligned}
\end{equation*}
where $\bm{w}(\xi,t)$ is the transverse displacement of the beam at position $\xi$ and time $t$, $u(t)$ is the acting force at time $t$, and $b(\cdot)$ is a given input distribution function describing how the control force acts along the spatial domain of the beam (e.g., $b(\xi)=\delta(\xi-\xi_0)$ models a point force at $\xi_0\in (0,\ell]$, see Figure~\ref{fig:beam_bnd}, while a smooth $b$ corresponds to a spatially distributed actuation), see Figure~\ref{fig:beam_distr}.
Further, $y(t)$ is the velocity that is co-located to the input force. The material parameters
$\rho,EI\in\Lp{\infty}([0,\ell])$ respectively stand for linear density and flexural rigidity. We assume that these functions are positive-valued with, moreover, $\rho^{-1},(EI)^{-1}\in\Lp{\infty}([0,\ell])$. Note that representing the Dirac $\delta$-distribution as an ordinary function is purely formal and does not meet the standards of full mathematical rigor. A more precise interpretation of such a configuration will be provided later.

The clamping condition at the left end means that 
\[{\bm w}(t,0)=\tfrac{\partial {\bm w}}{\partial\xi}(t,0)=0.\]
The other end is free, 
which means that both the bending moment and the shear force vanish, i.e.,
\[
  0=EI(\xi)\tfrac{\partial^2 \bm{w}}{\partial \xi^2}(\xi,t)\Big|_{\xi=\ell} =\tfrac{\partial}{\partial \xi}\!\left(EI(\xi)\tfrac{\partial^2 \bm{w}}{\partial \xi^2}(\xi,t)\right)\Big|_{\xi=\ell}.
\]
Our forthcoming representations use the {\em curvature} 
$\bm{\kappa}=\tfrac{\partial^2}{\partial\xi^2}\bm{w}$, and the 
{\em momentum density}
$\bm{p}=\rho\tfrac{\partial}{\partial t}\bm{w}$.

We discuss two scenarios: in the first, the beam is subject to a spatially distributed force with a co-located distributed velocity output; 
in the second, the beam is actuated by a point force at a specified location, and the output is given by the velocity of the displacement at that point. The system is then (still formally) represented by
{\allowdisplaybreaks\begin{align*}
\begin{pmatrix}\dot{\bm{\kappa}}(t,\xi)\\\dot{\bm{p}}(t,\xi)\end{pmatrix}
\!&=\!\!
\begin{bmatrix}0&\tfrac{\partial^2}{\partial\xi^2}
\\ -\tfrac{\partial^2}{\partial\xi^2}&0\end{bmatrix}\!\!\begin{pmatrix}EI(\xi)\bm{\kappa}(t,\xi)\\\rho(\xi)^{-1}\bm{p}(t,\xi)\end{pmatrix}\!+\!\begin{bmatrix}
0\\    
b(\xi)\end{bmatrix}u(t),\\
0&=\rho(\xi)^{-1}\bm{p}(t,\xi)\Big|_{\xi=0}=\tfrac{\partial}{\partial \xi}\!\left(\rho(\xi)^{-1}\bm{p}(t,\xi)\right)\Big|_{\xi=0},\\
0&=EI(\xi)\bm{\kappa}(\xi,t)\Big|_{\xi=\ell} =\tfrac{\partial}{\partial \xi}\!\left(EI(\xi)\bm{\kappa}(\xi,t)\right)\Big|_{\xi=\ell},\\
y(t)&=\int_0^\ell b(\xi)\rho(\xi)^{-1}\bm{p}(t,\xi)\dx[\xi].
\end{align*}}
We express the system in the framework of system nodes.
Our state space is $\Lp{2}([0,\ell];\C^2)$, equipped with the norm
\[\norm*{\pvek{\\[-7mm]\bm{\kappa}\\[-8mm]}{\bm{p}}}_\X=\left(\int_0^\ell EI(\xi)|\bm{\kappa}(\xi)|^2+\rho(\xi)^{-1}|\bm{p}(\xi)|^2\d\xi\right)^{1/2},\]
which is equivalent to the standard norm on $\Lp{2}([0,\ell];\C^2)$ due to our boundedness and strict positivity of $\rho$ and $EI$.\\
The system is single-input single-output, i.e., $\U = \Y = \C$.
To set up the system node, we
introduce the spaces
\begin{align*}
    \Hk{2}_{0l}([0,\ell])&=\setdef{f\in\Hk{2}([0,\ell])}{f(0)=f'(0)=0},\\
    \Hk{2}_{0r}([0,\ell])&=\setdef{f\in\Hk{2}([0,\ell])}{f(\ell)=f'(\ell)=0}.
\end{align*}
These will be used to formulate the main operator $A$, which is defined by 
$A:\X\supset\dom(A)\to\X$ with
\begin{align*}
\dom(A)&=\setdef*{\pvek{\\[-7mm]\bm{\kappa}\\[-8mm]}{\bm{p}}}{\;  \parbox{28mm}{$\rho^{-1}\bm{p}\in \Hk{2}_{0l}([0,\ell])$,\\ $EI\bm{\kappa}\in \Hk{2}_{0r}([0,\ell])$}},\\
A\pvek{\\[-7mm]\bm{\kappa}\\[-8mm]}{\bm{p}}&=\pvek{\phantom{-}\tfrac{\partial^2}{\partial\xi^2}\big(\rho^{-1}\bm{p}\big)}{-\tfrac{\partial^2}{\partial\xi^2}\big(EI\bm{\kappa}\big)}.
\end{align*}
It is straightforward to verify that $A$ is skew-adjoint, and hence, by \cite[Thm.~3.8.6]{TucsnakWeiss2009}, $A$ generates a strongly continuous semigroup (in fact, even a~group).
To formulate the input, we introduce the spaces
\begin{align*}
    \Hk{-2}_{0r}([0,\ell])&=\Hk{2}_{0l}([0,\ell])^*,\quad\Hk{-2}_{0l}([0,\ell])=\Hk{2}_{0r}([0,\ell])^*.
\end{align*}
It follows from the differentiation by parts formula that the second derivative on $\Hk{2}_{0l}([0,\ell])$ and $\Hk{2}_{0r}([0,\ell])$ respectively extend to bounded operators
\begin{align*}
    \tfrac{\partial^2}{\partial\xi^2}_l&\in\Lb(\Lp{2}([0,\ell]),\Hk{-2}_{0l}([0,\ell])),\\
    \tfrac{\partial^2}{\partial\xi^2}_r&\in\Lb(\Lp{2}([0,\ell]),\Hk{-2}_{0r}([0,\ell])).
\end{align*}
The extension of $A$ to $\X$ (cf.~Remark~\ref{rem:nodes}\,\ref{rem:nodes1}) is now given by 
\[A_{-1}=\begin{bmatrix}0&\tfrac{\partial^2}{\partial\xi^2}_r
\\ -\tfrac{\partial^2}{\partial\xi^2}_l&0\end{bmatrix}\begin{bmatrix}\rho^{-1}&0
\\ 0&EI\end{bmatrix}.\]
The input is modelled by a~distribution $b\in \H_{0r}^{-2}([0,\ell])$,
and, consequently, the input operator is given by
\[B=\bvek{0}{b}.\]
We now construct the system node in the opposite way to Remark~\ref{rem:nodes}\,\ref{rem:nodes1}, namely via 
\[A\&B\left(\begin{smallmatrix}
    \bm{\kappa}\\\bm{p}\\u
\end{smallmatrix}\right)=A_{-1}\left(\begin{smallmatrix}
    \bm{\kappa}\\\bm{p}
\end{smallmatrix}\right)+Bu,\]
and the domain of $A\&B$ is given by \eqref{eq:ABform}. The output operator
is defined by
\[C\&D\left(\begin{smallmatrix}
    \bm{\kappa}\\\bm{p}\\u
\end{smallmatrix}\right)=\scprod*{\rho^{-1}\bm{p}}{b}_{\Hk{2}_{0l},\Hk{-2}_{0r}},\]
where the latter denotes the duality pairing. From the above construction it follows that \eqref{eq:KYP} holds, and hence the system associated with the previously introduced system node is impedance passive.
For the choice of the force density distribution $b$, two cases are of particular interest:
\begin{enumerate}[label=(\alph{*})]
  \item $b \in \Lp{2}([0,\ell])$, corresponding to a {\em distributed control input}. In this case, the output is the weighted velocity average
  \[
    y(t) = \int_0^{\ell} \rho(\xi)^{-1}\bm{p}(t,\xi)b(\xi)\,d\xi,
  \]
where we tacitly assume that $b$ is real-valued, so as to avoid the need for complex conjugation. The domain of $A\&B$ is, in this case, given by $\dom(A)\times\C$.
\item $b = \delta_\xi$, where $\delta_\xi$ denotes the delta distribution at $\xi \in (0,\ell]$. This models a {\em point force applied at~$\xi$}.
Indeed, applying formula \eqref{eq:ABform} shows that the point force induces a jump in the derivative of $EI\bm{\kappa}$ at $\xi_0$. More precisely,
\begin{multline*}
\bigl[\tfrac{\partial}{\partial\xi}(EI\,\bm{\kappa})\bigr]_{\xi_0^-}^{\xi_0^+}\\
:= \tfrac{\partial}{\partial\xi}(EI\,\bm{\kappa})(\xi_0^+) - \tfrac{\partial}{\partial\xi}(EI\,\bm{\kappa})(\xi_0^-)
= u(t),
\end{multline*}
where the evaluation at $\xi_0^+$ and $\xi_0^-$ is to be understood as the limit taken from the right and from the left, resp. Even more precisely, the domain of $A\&B$ is given by
\begin{align*}
\dom(A\&B)\!&=\setdef*{\!\left(\begin{matrix}
    \bm{\kappa}\\\bm{p}\\u
\end{matrix}\right)\!}{\parbox{36mm}{$\rho^{-1}\bm{p}\in \Hk{2}_{0l}([0,\ell])$, $EI\bm{\kappa}\big|_{[0,\xi_0]}\!\in\! \Hk{2}([0,\xi_0])$, $EI\bm{\kappa}\big|_{[\xi_0,\ell]}\!\in\!\Hk{2}_{0r}([\xi_0,\ell])$, $\bigl[\tfrac{\partial}{\partial\xi}(EI\,\bm{\kappa})\bigr]_{\xi_0^-}^{\xi_0^+}=u$}},\\
A\&B\left(\begin{smallmatrix}
    \bm{\kappa}\\\bm{p}\\u
\end{smallmatrix}\right)&=\pvek{\phantom{-}\tfrac{\partial^2}{\partial\xi^2}\big(\rho^{-1}\bm{p}\big)}{-\tfrac{\partial^2}{\partial\xi^2}\big(EI\bm{\kappa}\big)},
\end{align*}
where the second derivatives above are to be under\-stood in such a way that any possible jump of $(EI\,\bm\kappa)_\xi$ at $\xi_0$ is ignored.
The output operator is given by the application of $\delta_{\xi_0}$ at $\rho^{-1}\bm{p}$, that is, the evaluation of $\rho^{-1}\bm{p}$ at $\xi_0$. In other words, the output consists of the displacement velocity at $\xi_0$.
Finally, we note that the choice $\xi_0=\ell$ leads to a boundary condition for $\tfrac{\partial}{\partial \xi}(EI\,\bm{\kappa})$ at $\ell$, and the corresponding velocity output is likewise located at the right end of the beam.\end{enumerate}

\subsection{Numerical simulations}

We set $\ell=1$ and take constant coefficients $EI\equiv \rho \equiv 1$. The reference signal is
\[
  y_{\mathrm{ref}}(t)=\cos t,
\]
and the funnel function is
\[
  \varphi(t)=10-9.5\,e^{-t/2},
\]
so that $1/\varphi(0)=2$ (funnel diameter $4$ at $t=0$) and $\lim_{t\to\infty}1/\varphi(t)=0.1$ (diameter $0.2$). We assume the beam is at rest at $t=0$, i.e., $x_0=0$, and that no external disturbances act on the system.
Two scenarios are considered:
\begin{enumerate}[label=(\alph{*})]
  \item\label{item:distr} \emph{Distributed actuation:} the input operator is given by the characteristic function on $[1/3,2/3]$, i.e., $b= 1_{[1/3,2/3]}$, with the output given by the spatial average velocity on $[1/3,2/3]$. By $x_0=0$ and $B\in\Lb(\U,\X)$, condition \eqref{eq:OutputInit} holds anyway, whence the 
  standard funnel controller
\[u(t)=-\frac{e(t)}{1-\varphi(t)^2\|e(t)\|^2}
\]
can be applied (no feedforward term according to Remark~\ref{rem:funnelass}\,\ref{rem:funnelass4} is required).
  \item\label{item:point} \emph{Point actuation:} a point force is applied at $\xi_0=\tfrac{1}{2}$, and the output is given by the displacement velocity at $\xi_0$.  
The initial output error satisfies \[e(0)=\rho(1/2)^{-1}\bm{p}(0,1/2)-\cos(0)=-1.\]
To compensate the mismatch at $t=0$, we equip our controller with a~feedforward term according to Remark~\ref{rem:funnelass}\,\ref{rem:funnelass4}. That is, our controller takes the form
\begin{equation}u(t)=\underbrace{-\frac{e(t)}{1-\varphi(t)^2\|e(t)\|^2}}_{=:u_{\rm fun}(t)}
+ \underbrace{p(t)\frac{e(0)}{1-\varphi(0)^2\|e(0)\|^2}}_{=:u_{\ext}(t)}\label{eq:uEBBcont}
\end{equation}
for
\[p(t)=\begin{cases}\tfrac12\big(1+\cos(\pi t)\big)&:\;t\in[0,1],\\0&:\;t\in(1,\infty).\end{cases}\]
Note that, by using $e(0)=-1$ and $\varphi(0)=1/2$,
\[u_{\ext}(t)=\begin{cases}-\tfrac23\big(1+\cos(\pi t)\big)&:\;t\in[0,1],\\0&:\;t\in(1,\infty).\end{cases}\]

\end{enumerate}

\begin{figure}[ht]
  \centering
  \input{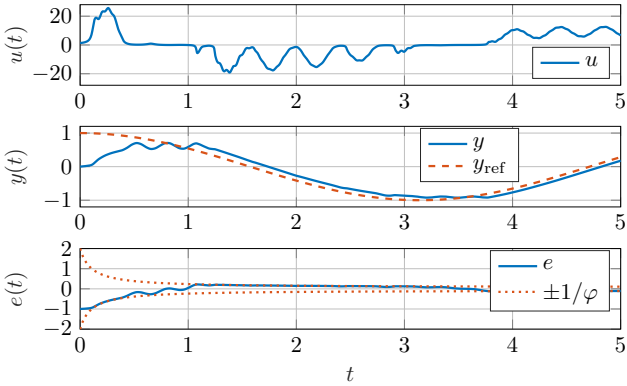}
  \caption{Scenario \ref{item:distr}: Distributed actuation and observation}\label{fig:distr}
\end{figure}

\begin{figure}[hb]
  \centering
  \input{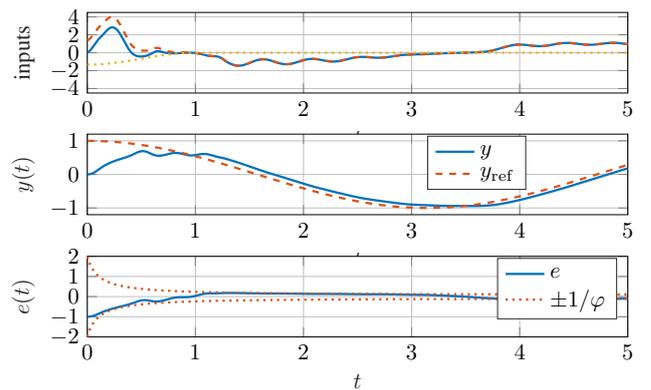}
  \caption{Scenario \ref{item:point}: Point actuation and observation, with input $u=u_{\rm fun}+u_{\rm ext}$ ($u_{\rm fun}$ dashed red, $u_{\rm ext}$ dashed yellow, $u$ solid blue), with notation as in \eqref{eq:uEBBcont}}\label{fig:point}
\end{figure}

In both cases, Theorem~\ref{thm:Main_Thm_Exist_Uniq} guarantees that, under the funnel controller, the tracking error $e(t)=y(t)-y_{\mathrm{ref}}(t)$ remains within the funnel, i.e., $e(t)\in(-1/\varphi(t),1/\varphi(t))$. For our simulations (see Figure~\ref{fig:distr} and Figure~\ref{fig:point}) we discretize the beam using cubic Hermite ($\Hk{2}$) finite elements on an equidistant mesh with $k=80$ nodes, so that all elements have uniform length $1/(n-1)$. This results in a~semidiscretized system with state space dimension $n=324$.

\balance
\bibliographystyle{abbrv}        
\bibliography{References}           

\end{document}